\documentclass[11pt]{article}
\usepackage{amssymb}
\usepackage{latexsym}
\newtheorem{thm}{Theorem}[section]
\newtheorem{prop}[thm]{Proposition}
\newtheorem{lemma}[thm]{Lemma}
\newtheorem{cor}[thm]{Corollary}

\newcommand{\call}{{\cal L}}

\newcommand{\cinf}{C^{\infty}}

\def\qed{\rule{2.3mm}{2.3mm}}

\newcommand{\lag}{\langle}
\newcommand{\rg}{\rangle}

\begin{document}

\title{\bf On the local structure of Dirac manifolds}

\author{
 {\bf Jean-Paul Dufour} \thanks{email: dufourj@math.univ-montp2.fr}\\
 D\'epartement de Math{\'{e}}matiques \\
 CNRS-UMR 5030, Universit\'e Montpellier 2\\
 Place Eug\`ene Bataillon, 34095 Montpellier Cedex 05\\
\\[.5cm]
 {\bf  A\"{\i}ssa Wade} \thanks{email: wade@math.psu.edu}
 \thanks{Research partially supported by the Shapiro Funds} \\
 Department of Mathematics\\
  The Pennsylvania State University \\
 University Park, PA 16802 \\
}

\date{}
\maketitle

\begin{abstract}
 We give a local normal form for Dirac structures.
 As a consequence, we show that the
 dimensions of the pre-symplectic leaves
 of a Dirac manifold have the same parity.
 We also show that, given a point $m$ of a Dirac manifold $M$,
 there is a well-defined
 transverse Poisson structure to the pre-symplectic leaf $P$
 through $m$. Finally, we describe the neighborhood of a pre-symplectic leaf
 in terms of geometric data. This description agrees with that given by
 Vorobjev for the Poisson case.
\end{abstract}

\section{Introduction}
 A {\em Dirac manifold}
  is a smooth manifold $M$ equipped with a vector subbundle
$L$ of the Whitney sum $TM \oplus T^*M$  which is 
 maximal isotropic with respect to the natural pairing on $TM \oplus T^*M$
 and integrable in the sense that the smooth sections 
 of $L$ are closed under the Courant bracket 
 (see Section 2).
 The vector bundle $L$ is then called a {\em Dirac structure}
 on $M$. 

\medskip
 Dirac structures on manifold were first introduced by Courant 
 and Weinstein in the mid-eighties \cite{CW86}.
 A few years later, further investigations
 were undertaken  in \cite{C90}. Recently, the theory of Dirac structures has been extensively developed  in connection with 
 various topics in mathematics and physics (for instance, see 
\cite{BC97}, \cite{G04}, \cite{LWX97},  \cite {BW04}). Specific examples 
 of Dirac manifolds include pre-symplectic and Poisson manifolds.
 Thus, it is important to understand the local structure of
 a Dirac manifold. The main goal of this paper is to provide a description
  of the local structure of such a manifold.

 Every Dirac manifold admits a foliation by pre-symplectic
 leaves. To our knowledge, the local structure of Dirac manifolds
 was studied only  in neighborhoods of
  regular points \cite{C90} (see also \cite{G04} for
 the case of complex Dirac structures).
  By a regular point we mean a point  for which there is an
 open neighborhood where the foliation is regular.
 It is natural to ask  about the local structure around
  non regular points.
 In Section 3, we
 give a normal form for a Dirac structure $L$ on a smooth manifold
 $M$  near an arbitrary point $m \in M$
 (see Theorem \ref{normal form}). This normal form allows us to
 conclude that the dimensions of the pre-symplectic leaves 
 have the same parity.  

\medskip

 We show in Section 4 that,
  given a point $m$ in  a Dirac manifold $M$,  there  is a well-defined
 transverse Poisson structure whose rank at $m$ is zero.
 This extends facts from the classical case of Poisson structures 
 (see \cite{We83}).  In Section 5, we describe the neighborhood
 of a pre-symplectic leaf of a Dirac manifold using
 the concept of a geometric data (see \cite{V00}).
  Dirac structures on manifolds are constructed from
  given geometric data. We prove that, conversely, 
 one can construct a geometric data from a Dirac manifold
 $M$ with a fixed tubular neighborhood of a symplectic
 leaf $P$.

\medskip

 The paper is divided into five sections.
 Section 2 contains some basic definitions and results.
  Our  main theorems are in Sections 3-5
 (see Theorems \ref{normal form}, \ref{transversal},
 \ref{geometric data}, and \ref{Part 2}).

\section{Preliminaries}

 Let $M$ be a smooth $n$-dimensional
manifold. We denote by $ \lag  \cdot , \cdot \rg$ the canonical
 symmetric bilinear operation
  on the vector bundle $TM \oplus T^*M \rightarrow M$.
 This induces a symmetric $\cinf$-bilinear
 operation on the space of sections of $TM \oplus T^*M$ given by:
$$ \lag (X_1, \alpha_1), (X_2, \alpha_2) \rg 
={1 \over 2}(i_{X_2} \alpha_1 + i_{X_1} \alpha_2), \ \forall \
 (X_1, \alpha_1),  \ (X_2, \alpha_2) \in
 \Gamma(TM \oplus T^*M).$$

 An {\em almost Dirac structure} on $M$ is a subbundle of
$TM \oplus T^*M \rightarrow M$ which is  maximal isotropic
 with respect to $ \lag  \cdot , \cdot \rg$.
\medskip

\noindent The {\em Courant bracket}
 on  $\Gamma(TM \oplus T^*M)$ is defined by:
$$[(X_1, \alpha_1), \ (X_2, \alpha_2)]_{_C}=([X_1,X_2], \
  \call_{X_1}\alpha_2 - i_{X_2}d\alpha_1),$$
\noindent  where $\call_X=d \circ i_X + i_X \circ d$
  is the Lie derivation by $X$.

 A {\em Dirac structure} $L$ on $M$
 is an almost Dirac structure which is integrable (i.e.
  $ \Gamma(L)$  is closed under the Courant bracket). In this case,
   the pair $(M, L)$  is called a
 {\em  Dirac manifold}.

\bigskip

\noindent{\bf Examples}

\noindent {\bf (i)} Let $\Omega$ be a 2-form on $M$.
 Consider the graph $$L_{\Omega}=\{ (X, i_X \Omega) \ | \  X \in TM\}.$$
Then $L_{\Omega}$ is a Dirac structure if and only if
$d{\Omega}=0.$ Furthermore,  a Dirac structure is the graph
 of a 2-form if and only if  $L \cap (\{0\} \oplus T^*M)=\{0\}$
 at every point.

\medskip

\noindent {\bf (ii)} Let $\pi$ be a bivector field on $M$. We denote
$$L_{\pi}= \{ (\pi \alpha, \alpha), \ | \ \alpha \in T^*M \}.$$
\noindent Then $L_{\pi}$ is Dirac if and only if
 $\pi$ is a Poisson tensor.
 Furthermore,  a Dirac structure is the graph
 of a bivector field if and only if  $L \cap (TM\oplus \{0\})=\{0\}$
 at every point.

\hfill \qed

\bigskip

  Let $L$ be an almost Dirac structure  on $M$.
  Consider the distribution
 $$({\cal D}_L)_x= pr_1(L_x) \quad \mbox{for all} \ x \in M,$$
\noindent  where $pr_1$ is the canonical projection of
 $L_x$ onto $T_xM$.
 The distribution ${\cal D}_L$ is involutive
 when $L$ is integrable.
 Hence, in this case, $L$ gives rise to a  singular foliation.
 Furthermore,  there is  a skew-symmetric
 bilinear map
$\Omega_{_L}:  {\cal D}_L \times {\cal D}_L \rightarrow
 \cinf(M)$ given by
 $$\Omega_{_L}(X,Y)= \alpha(Y), \ \ \mbox{\rm for any} \
 (X, \alpha), (Y, \beta) \in \Gamma(L). \leqno(1)$$

\noindent We have the following proposition:

\begin{prop}{\rm (~\cite{C90})}
\label{prop1}
 If $L$ is a  Dirac structure on $M$
 then $d\Omega_{_L}=0$.
\end{prop}

 \noindent {\bf Remark.} As an immediate consequence of this proposition,
  one sees that
 a  Dirac structure on $M$ gives rise to a singular foliation
 by pre-symplectic leaves (i.e.  on each leaf, there is a closed 2-form).

\section{ A local normal form for Dirac manifolds}

\begin{prop}
\label{Poisson}
Let $L$ be a Dirac structure a on smooth manifold $M$
 of dimension $n$, and let $m_0 \in M$.
 If the pre-symplectic leaf  through $m_0$ is a single point
  then there is a neighborhood
 $U$ of $m_0$ such that $L_{|U}$ is the graph of a
  Poisson structure $\Pi$.
\end{prop}

\noindent {Proof:}
Assume that the pre-symplectic leaf  through $m_0$ is $P=\{m_0\}$.
There are vector fields $X_1 , \dots , X_n$, and
 1-forms $\alpha_1, \dots , \alpha_n$ defined in an
  open neighborhood $U$  of $m_0$ such that $L_{|U}$
 is determined by the local sections $e_i=(X_i, \alpha_i)$,
 where $X_i(m_0)=0$. In local coordinates
 $(y_1, \dots, y_n)$ such that $y_i(m_0)=0$, we have
 the expressions
$$X_i= \sum_{i=1}^n X_{ij}(y) {\partial \over \partial y_j},
 \quad  \alpha_i=\sum_{i=1}^n \alpha_{ij} dy_j.$$
 The $e_i(m)$ give the following matrix
 $$ \left( \begin{array}{cccccccc}
  X_{11} & \dots &X_{1n} & \alpha_{11} & \dots & \alpha_{1n} \cr
 \vdots & & & & & \cr
 X_{n1} & \dots &X_{nn} & \alpha_{n1} & \dots & \alpha_{nn}
\end{array}\right) \quad \mbox{ with $X_{ij}(m_0)=0$}.$$
\noindent  The sub-matrix
 $(\alpha_{ij}(m_0))$ is invertible since  $dim (L_{m_0})=n$.
 Therefore, $(\alpha_{ij}(m))$ remains invertible at all points
  in a small neighborhood of $m_0$. Let
 $(\alpha^{ij}(m))$ be the inverse of $(\alpha_{ij}(m))$.
 Define $$e'_i= \sum_{j=1}^n \alpha^{ij} e_j, \quad \forall i= 1 , \dots , n.$$
This can be written as
$$e'_i=\Big(\sum_{i=1}^n X'_{ij}{\partial \over \partial y_j},  \ dy_i\Big)
  \quad \mbox{with} \quad  X'_{ij}=-X'_{ji}.$$
\noindent Define
$$\Pi= \sum_{i<j} X'_{ij}{\partial \over \partial y_i} \wedge
{\partial \over \partial y_j}.$$
\noindent Then, the Schouten bracket $[\Pi, \Pi]$ is zero
  (this is due to the fact that
 $\Gamma(L)$ is closed under the Courant bracket).
Furthermore, $L_{|U}$  is the graph of $\Pi$.

 \hfill \qed

 Now, we assume that the  pre-symplectic leaf
 $P$   through $m_0 \in M$ is not a single point.
  Let $U$ be a  neighborhood of $m_0$
  with local coordinates
 $(x_1, \dots, x_r, y_1, \dots , y_s)$
   such that $x_i(m_0)=y_j(m_0)=0$,
 for all
  $i \in \{1, \dots, r\}$, $j\in \{1, \dots , s\}$, and such that the pre-symplectic leaf
  through $m_0$ has equations $y_1=0,\dots ,y_s=0.$ In the sequel, we will use the following notations
  $$x=(x_1, \dots, x_r),\ \ y=( y_1, \dots , y_s).$$
  
  \noindent Without loss of generality, we can assume that
 there are vector
  fields $Y_i(x,y)$,
 $Z_j(x,y)$ defined on $U$ and 1-forms  $\alpha_i(x,y)$, $\beta_j(x,y)$  such that
 $\Gamma(L_{|U})$ is spanned by $$ {\cal S}_i=\Big(
 {\partial \over \partial x_i}+ Y_i(x,y), \ \alpha_i(x,y)\Big),
  \quad {\cal T}_j=(Z_j(x,y), \ \beta_j(x,y)),$$
 \noindent  with $Y_i(x,0)=0$, $Z_j(x,0)=0$, for all
  $i=1, \dots, r$, $j=1,  \dots , s$.
 We want to find  a new spanning set of
  local sections $\{ {\cal S}'_i, {\cal T}'_j \}$
  defined around $m_0$ which have very simple  expressions.
 In other words, we want  to find a normal form
of $L$  at $m_0$.
   Denote
$$ X_i= {\partial \over \partial x_i}+ Y_i(x,y)=
 {\partial \over \partial x_i}+  \sum_{j=1}^r {\widehat Y}_{ij} (x,y)
 {\partial \over \partial x_j}+ \sum_{j=1}^s{\widetilde Y}_{ij} (x,y)
  {\partial \over \partial y_j}.$$

\noindent  One can notice that there are smooth functions
 $f_{ij}(x,y)$ such that
$$\sum f_{ij}(x,y) X_j= {\partial \over \partial x_i}
+ \sum_{j=1}^sY'_{ij} (x,y)
  {\partial \over \partial y_j}.$$
 Equivalently, there is a  matrix $(f_{ij}(x,y))$
 whose coefficients are smooth functions
  and such that

$$ \left( \begin{array}{c}
  I+({\widehat Y}_{ij})
\end{array}\right)
 \left( \begin{array}{c}
 f_{ij}\end{array}\right) =I
$$

\noindent Indeed,
 $( I+({\widehat Y}_{ij}))$
  is invertible at all
 points in $U$ (up to a shrinking of $U$). Hence, there  are smooth functions
 $f_{ij}(x,y)$ satisfying the above matrix equation.
 It follows that  $\Gamma(L_{|U})$ is spanned by smooth sections
 of the form  $${\cal S}'_i=\Big({\partial \over \partial x_i}
+ \sum_{j=1}^sY'_{ij} (x,y)
  {\partial \over \partial y_j}, \  \alpha'_i(x,y)\Big), \quad
 {\cal T}_j=\Big(Z_j(x,y), \ \beta_j(x,y)\Big),$$
\noindent for $i=1, \dots, r,$ and $j=1, \dots, s.$
 We write
$$Z_i=
 \sum_{j=1}^r {\widehat Z}_{ij} (x,y)
 {\partial \over \partial x_j}+ \sum_{j=1}^s{\widetilde Z}_{ij} (x,y)
  {\partial \over \partial y_j}.$$
 Define
$${\cal T}'_i=  {\cal T}_i -
 \sum_{j=1}^r {\widehat Z}_{ij}(x,y)  {\cal S}'_j.$$
 Then, we see  that
 $\Gamma(L_{|U})$ is spanned by  smooth sections
 of the form  $${\cal S}'_i=\Big({\partial \over \partial x_i}
+ \sum_{j=1}^s Y'_{ij}(x,y)
  {\partial \over \partial y_j}, \  \  \alpha'_i(x,y)\Big), \quad
{\cal T}'_j=\Big(\sum_{k=1}^s Z'_{jk}(x,y)
 {\partial \over \partial y_k}, \  \ \beta'_j(x,y)\Big), $$

\noindent where  $Y'_i(x,0)= Z_j (x, 0)=0$, for
 all $i=1, \dots, r, \quad j=1, \dots, s.$
 Using the fact that $L$ is isotropic, we get
 $$ \alpha'_i ({\partial \over \partial x_j})+
  \alpha'_j ({\partial \over \partial x_i}) =0 \quad \mbox{and} \quad
 \beta'_j( {\partial \over \partial x_i})=0
\quad \mbox{at every point } \ p  \in P.$$
\noindent Moreover, a basis for
 the fiber $L_{m_0}$ is given by the elements

$${\cal S}'_i (m_0)= ({\partial \over \partial x_i}, \
 \alpha'_i), \quad
 \quad {\cal T}'_j(m_0)=\Big(0, \ \sum_{k=1}^s \beta'_{jk}
 {\partial \over \partial y_k}\Big),$$
\noindent for $i=1, \dots, r, \quad j=1, \dots, s.$
 Using matrix notations, we  can put the ${\cal S}'_i(m_0)$ and
 ${\cal T}'_j(m_0)$ into row vectors which give
  the following rectangular matrix:

$$ \left( \begin{array}{cccc}
  I& 0 & *&
  *  \cr
  {}\cr
  0 & 0 & 0& (\beta'_{ij})
\end{array}\right).
$$
\noindent The sub-matrix $(\beta'_{ij}(m_0))$ is invertible
 since $dim (L_{m_0})=r +s=dim M$. Hence,  $(\beta'_{ij}(x,y))$
  is invertible at all points $m=(x,y)$ in a small neighborhood of $m_0$.
 Let $(g_{ij} (x,y))$ be the inverse of the matrix
 $(\beta'_{ij}(x,y))$. We denote
$${\cal T}''_i= \sum_{j=1}^s g_{ij}(x,y){\cal T}'_j.$$
 Hence, ${\cal T}''_i$ has the form
 $$ {\cal T}''_i=\Big(\sum_{k=1}^s Z''_{ik}(x,y)
 {\partial \over \partial y_k}, \  dy_i+
 \sum_{k=1}^r \beta''_{ik}(x,y) d x_k \Big).$$
 \noindent Now, we replace
 $${\cal S}'_i=\Big({\partial \over \partial x_i}
+ \sum_{j=1}^sY'_{ij} (x,y)
  {\partial \over \partial y_j}, \  \ \sum_{j=1}^r \overline{\alpha}'_{ij}
 (x,y) dx_j + \sum_{k=1}^s \widetilde{\alpha}'_{ik}(x,y) dy_k \Big)$$
 \noindent by the following
 $${\cal S}''_i={\cal S}'_i - \sum_{k=1}^s \widetilde{\alpha}'_{ik}(x,y)
 {\cal T}''_k.$$

\noindent We then obtain a new spanning  set
  $\{ {\cal S}''_i, {\cal T}''_j\}$ of local
 sections of $L_{|U}$.
  We summarize the above discussion in the following theorem:
\begin{thm}
\label{normal form}
Let $L$ be a Dirac structure on a smooth manifold $M$.
 Given any point $m_0 \in M$, there is a coordinate system
 $(x_1, \dots , x_r, y_1 , \dots , y_s)$ defined on
 an open neighborhood $U$ of $m_0$ such that the intersection
 of $U$ with the pre-symplectic leaf through
 $m_0$ is  the set $\{y_1= \dots = y_s=0\}$, and
 $\Gamma(L_{|U})$ is spanned by sections of the form
$$
 {\cal H}_i=\Big({\partial \over \partial x_i}
+ \sum_{k=1}^sX_{ik} (x,y)
  {\partial \over \partial y_k}, \ \   \sum_{k=1}^r\alpha_{ik}(x,y) dx_k
 \Big)$$
\noindent and
$${\cal V}_j=\Big(\sum_{k=1}^s Z_{jk}(x,y)
 {\partial \over \partial y_k}, \  \ dy_j+
 \sum_{k=1}^r \beta_{jk}(x,y) d x_k \Big),$$
where $X_{ik}(m_0)=0$,  $ Z_{jk}(m_0)=0$, for all
 $i \in \{1 , \dots , r \}$ and  \ $j \in \{1 , \dots , s \}$.
\end{thm}

We will use the following notations
$$X_i={\partial \over \partial x_i}
+ \sum_{k=1}^sX_{ik} (x,y)
  {\partial \over \partial y_k},\ \quad \alpha_i =\sum_{k=1}^r\alpha_{ik}(x,y) dx_k,$$
  $$\ \ Z_j=\sum_{k=1}^s Z_{jk}(x,y)
 {\partial \over \partial y_k},\ \quad \beta_j=\ dy_j+
 \sum_{k=1}^r \beta_{jk}(x,y) d x_k .$$

\noindent{\bf Remark.}

\noindent {\bf a)} The normal form
 $ ({\cal H}_i, {\cal V}_j)$ persists  when
 a change of coordinates of the type
 $\Phi(x,y)= (x, \Phi_2(x,y))$
 is performed.

\medskip

\noindent {\bf  b)} Using the  fact that $L$ is isotropic, one gets
 $$ 2 \langle {\cal H}_i, \ {\cal V}_j \rangle= X_{ij}+ \beta_{ji}=0,$$
\noindent  for all  $i \in \{1 , \dots , r \}$ and
  \ $j \in \{1 , \dots , s \}$. Hence $X_{ij}=- \beta_{ji}$.
Furthermore,
$$ 2 \langle {\cal V}_i, \ {\cal V}_j \rangle=
 Z_{ij}+Z_{ji}=0.$$
\noindent This shows that the matrix $(Z_{ij}(x,y))$ is skew-symmetric.
An analogous calculation shows that the matrix $(\alpha_{ij}(x,y))$ is also skew-symmetric.

\medskip

\noindent {\bf  c)}
 If $\Omega_L$ is the 2-form
 associated with the Dirac  structure $L$
 (it is defined as in (1)) then we have
$$ \Omega_L\Big(X_i, \ X_j \Big)= \alpha_{ij},
 \quad \Omega_L\Big( X_i, \  Z_j\Big)=0, \quad
\Omega_L(Z_i, \ Z_j)=Z_{ij}.$$

\begin{cor}
Given a Dirac structure $L$ on a smooth manifold $M$, the
  dimensions of the leaves of its
 associated pre-symplectic foliation
 have the same parity.
\end{cor}

\noindent {\it Proof:}
  It is sufficient to work in a small
 neighborhood of an arbitrary point $m_0 \in M$.

\noindent {\it Case 1:} Suppose that the leaf through $m_0$ is a single point.
 By proposition \ref{Poisson}, $L$ is the graph of a Poisson structure
  in a  neighborhood of $m_0$. Hence, the dimensions of the pre-symplectic
 leaves are even.

\medskip
\noindent {\it Case 2:}
 Suppose that the leaf through $m_0$ is  not a single point.
 Then,  Theorem \ref{normal form} provides
   a  coordinate system $(x_1, \dots , x_r, y_1, \dots
 y_s)$  defined on an open neighborhood $U$
 of $m_0$ and  a spanning
 set $\{ {\cal H}_i, \ {\cal V}_j \}$ of local sections
 of $L_{|U}$ such that
 $$pr_1({\cal H}_i)=X_i={\partial \over \partial x_i}
+ \sum_{k=1}^sX_{ik} (x,y)
  {\partial \over \partial y_k}, \quad
 pr_1({\cal V}_j)=Z_j=\sum_{k=1}^s Z_{jk}(x,y)
 {\partial \over \partial y_k},$$
 \noindent where the functions $X_{ik}$, $Z_{ik}$
 vanish at $m_0$.
 Let $({\cal D}_L)_m$ be the tangent space
  to the leaf through $m=(x,y)$. It is spanned by the set of vectors
 $\{X_i(x,y), Z_j(x,y) \ | \ i=1 , \dots , r \ j=1 , \dots ,  s\}$
  which corresponds to
 the matrix
 $$ \left( \begin{array}{ccc}
  I & & * \cr
  0 &  & (Z_{jk}(x,y))
\end{array}\right).$$
 The sub-matrix $(Z_{jk}(x,y)) $ is skew-symmetric, hence
 its rank is even. Since
 $$dim({\cal D}_L)_m= r + {\rm rank }(Z_{jk}(x,y)),$$
\noindent  we conclude that $dim({\cal D}_L)_m$
 and $dim({\cal D}_L)_{m_0}=r$ have the same parity.

 \hfill \qed

\medskip We should mention this 
  phenomenon of jumping dimensions
  appeared in the study of generalized
 complex structures on even-dimensional manifolds
 (see \cite{G04}). Moreover,
 as A. Weinstein points out (private communication),  this corollary
 can be obtained from Lemma 2.2
 of the paper \cite{T-W03}.

\section{Transverse Poisson structures}

\subsection{Induced Dirac structures on submanifolds}

Let $L$ be a Dirac structure on a manifold $M$.
 Let $Q$ be a submanifold of $M$.  In this section,
we will review a result  established in \cite{C90}
 which says that, under  certain regularity
 conditions, $L$ induces a Dirac structure on $Q$.
At every point $q \in Q$, we get a maximal isotropic vector space
 $$(L_Q)_q={ L_q \cap (T_qQ  \oplus T_q^*M) \over  L_q \cap
(\{0\} \oplus T_q Q^{\circ})},$$
\noindent where $T_q Q^{\circ}=\{ v \in T^*_qM \ | \
 v_{|T_qQ}=0 \}.$
\noindent Using the  map
 $(L_Q)_q \rightarrow T_qQ  \oplus T_q^*Q$ given by
$$(u,v) \mapsto (u, v_{|T_qQ}),$$
\noindent one can identify $(L_Q)_q$ with a subspace of
 $T_qQ  \oplus T_q^*Q$.
 In fact, $L_Q$ defines a smooth subbundle of
 $TQ \oplus T^*Q$ if and only if $L_q \cap (T_qQ  \oplus T_q^*M)$
  has constant dimension. Moreover, one has the following result:

\begin{prop}{\rm (see \cite{C90})}
\label{induced}
If $L_q \cap (T_qQ  \oplus T_q^*M)$
  has constant dimension  then $L_Q$ is a Dirac structure on $Q$.
\end{prop}
\subsection{Existence of a transverse Dirac structure}

 Let $L$ be a Dirac structure on $M$,
  $m_0$  a point of  $M$ and $Q$ a submanifold of $M$ which contains $m_0$
  and is transversal to the pre-symplectic leaf of $m_0$ in the sense that
   the tangent space of $M$ at $m_0$ is the direct sum of the tangent spaces of
  $Q$ and of the pre-symplectic leaf $P$.
  Choose  coordinates $(x_i,y_j)$ defined on
  an open neighborhood $U$ of $m_0$ as in Theorem
\ref{normal form} but with the additional condition that $Q$ is given by
 equations $x_1=0,\dots ,x_r=0.$ We adopt notations of the preceding section, i.e.
 $\Gamma(L_{|U})$ is spanned by sections of the form
$$
 {\cal H}_i=\Big({\partial \over \partial x_i}
+ \sum_{k=1}^sX_{ik} (x,y)
  {\partial \over \partial y_k}, \ \  \sum_{k=1}^r  \alpha_{ik}(x,y)
 d x_k \Big)$$
\noindent and
$${\cal V}_j=\Big(\sum_{k=1}^s Z_{jk}(x,y)
 {\partial \over \partial y_k}, \  \ dy_j+
 \sum_{k=1}^r \beta_{jk}(x,y) d x_k \Big),$$
where $X_{ik}(m_0)=0$,  $ Z_{jk}(m_0)=0$,
 for all $i \in \{1 , \dots , r \}$ and  \ $j \in \{1 , \dots , s \}$.

\begin{lemma}
 The  vector spaces  $L_q \cap (T_qQ \oplus T_q^*M)$  have the same dimension,
 for all $q \in Q$.
\end{lemma}

\noindent {\it Proof:}
Suppose that  $(u (q), \ v(q))$ is a vector in
$ L_q \cap (T_qQ \oplus T_q^*M).$ We write
$$(u (q), \ v(q))= \sum_{i=1}^r \lambda_i {\cal H}_i(q)+
\sum_{j=1}^s \mu_j {\cal V}_j(q).$$
 \noindent Then,
$$dx_k\Big(\sum_{i=1}^r \lambda_i({\partial \over \partial x_i}+ Y_i) +
\sum_{j=1}^s \mu_j Z_{jk} {\partial \over \partial y_k} \Big)=
 \lambda_k= 0.$$

\noindent Consequently,  $dim (L_q \cap (T_qQ \oplus T_q^*M)) \leq s$.
 But, the vectors $({\cal V} (q))_{ j=1, \dots s}$ are linearly independent
 at $q=m_0$. Therefore,
   they are linearly independent for all $q \in Q$
 (we can suppose that the open neighborhood
  $U$ of $m_0$  is  small enough).
 This shows that $L_q \cap (Vert_q \oplus T_q^*M)$ has constant dimension.

\hfill \qed

 Now, applying Proposition~\ref{induced},
 one can conclude that  $L_Q$ is
 a Dirac structure on $Q$.
 In fact, $L_Q$ is spanned by the sections
$${\cal V}_j(q)=\Big( Z_{j}
 , \  \ dy_j \Big)(q), \quad \forall q \in Q,
 \ \ j \in \{1 , \dots , s \},$$
 where we use here notations of the previous section.

Since the pre-symplectic leaf of $L_Q$ at $m_0$ reduces to a point, 
 the proposition \ref{Poisson}
 shows that $L_Q$ is the graph of  a Poisson structure.
 The corresponding Poisson tensor is given by
  $$\Pi_{Q}(dy_i, dy_j)= - Z_j(0,y) \cdot y_i.$$
  We have then proved the following result.

  \begin{thm} Let $Q$ be a submanifold  transversal to a pre-symplectic leaf
  $P$ of the Dirac manifold $M$ at a
  point $m_0$ ($T_{m_0}M=T_{m_0}P\oplus T_{m_0}Q$).
 Then the Dirac structure induces a Poisson structure $\Pi_Q$ on a
  neighborhood of $m_0$ in $Q,$ with $\Pi_Q(m_0)=0.$
  \end{thm}

  The above calculations show also that there is an induced
  Poisson structure on each submanifold
  given by equations $x=$constant. These  Poisson structures fit together
  to give a Poisson tensor $\Pi^V$ defined on a  whole neighborhood of $m_0$
  in $M$ by
$$\Pi^V(dy_i, dy_j)= - Z_j(x,y) \cdot y_i,$$
\noindent where $y=(y_1 , \dots , y_s)$ (respectively 
 $ x=(x_1, \dots , x_r)$) are local coordinates of $Q$ 
 (respectively $P$) around $m_0$.
\begin{lemma}
\label{prep 2}
 For any $i=1, \dots , r$, we have
$$[X_i, \ \Pi^{V}]=0,$$
\noindent where $X_i=p_1({\cal H}_i)$.
\end{lemma}

\noindent {\it Proof:}  Recall that  $ {\cal V}_j=(Z_j,\beta_j).$ For simplicity, we write $ \beta_j$ on the form
$dy_j+\beta^V_j.$ Then, we have
$$[{\cal H}_i, {\cal V}_j]=
 ([X_i, Z_j], \ d(X_i \cdot y_j) + \call_{X_i} \beta^V_j
 - i_{Z_j}d \alpha_i).$$

\noindent The fact that $L_{|U}$ is isotropic implies
\begin{eqnarray*}
0&=&2 \langle [{\cal H}_i, {\cal V}_j], \  {\cal V}_k \rangle \cr
 &=&[X_i, Z_j] \cdot y_k + Z_k \cdot (X_i \cdot y_j) +
 Z_k \cdot (\beta^V_j(X_i)) + d \beta^V_j (X_i, Z_k)
 - d \alpha_i(Z_j, Z_k) .
\end{eqnarray*}
\noindent But
$d \alpha_i(Z_j, Z_k) =0$ because
  $Z_j$ and $Z_k $  have only terms in $\frac{\partial}{\partial y}.$
 Moreover,
 if we denote
 $$\beta^V_j= \sum_{k=1}^r \beta_{jk}(x,y) dx_k,$$
\noindent then
 $$ Z_k \cdot (\beta^V_j(X_i)) + d \beta^V_j (X_i, Z_k)
  =
 \sum_{\ell =1}^s Z_k^\ell {\partial \beta_{ji} \over \partial
 y_{\ell}} -  \sum_{\ell =1}^s Z_k^\ell
 {\partial \beta_{ji} \over \partial y_{\ell}}=0.$$

\noindent There follows
$$2 \langle [{\cal H}_i, {\cal V}_j], \  {\cal V}_k \rangle =
 [X_i, Z_j] \cdot y_k + Z_k \cdot (X_i \cdot y_j)=0.$$

\noindent This equation can be written as
$$ X_i \cdot (Z_j \cdot y_k)
 - Z_j \cdot (X_i \cdot y_k)+ Z_k \cdot (X_i \cdot y_j)=0.$$

\noindent This is equivalent to the equation
 $$[X_i, \ \Pi^{V}] (dy_j, dy_k)=0,$$

 \noindent for all
 indexes $ i   \in \{1 , \dots , r\},$ and
 $j   \in \{1 , \dots , s\}.$
This completes the proof of the lemma.

\hfill \qed

\begin{thm}
\label{transversal}
 Let $Q$ and $Q'$ two submanifolds transversal to a same pre-symplectic
leaf $P$ of a Dirac structure $L.$ The Poisson structures induced by $L$ on $Q$  and $Q'$ 
 (near $Q\cap P$ and 
$Q'\cap P$ respectively) are locally isomorphic.
\end{thm}

\noindent {\it Proof:} By connexity of $P,$ it is sufficient to construct the isomorphism in the case where $Q$
and $Q'$ are near enough. We can also suppose that $Q\cap P$ and $Q'\cap P$ are different. 
Hence, it is enough to work in a domain with coordinates $(x,y)$ as above,
 that is,  $P$ has equation $y=0,$ $Q$ has equation $x=0$
and, moreover, $Q'$ has equation $x=x^0$ where $x^0$ is some constant different from $0.$ Now we will use
 Lemma \ref{prep 2}: because $X_i$ has a component $\frac{\partial}{\partial x_i}$ we can go from $0$ to $x^0$ in $P$
using a sequence of trajectories of the different fields $X_i$, moreover the flows of these fields preserve
verticals $x=$constant and Lemma \ref{prep 2} says that they preserve also $\Pi^V,$ so they exchange the Poisson structures
on the verticals. \hfill \qed

\medskip
\noindent{\bf Remark:}
 It follows from Theorem \ref{transversal} that each pre-symplectic leaf of a Dirac structure has a well defined, up  to isomorphism, Poisson transversal
structure. This extends a well-known result in the  Poisson case.
\medskip

We can also remark that, in the classical case of Poisson structures, the above method  used to prove the uniqueness of the
 transverse Poisson structure is simpler than the ones in literature 
 (see for instance \cite{We83}).

\section{Geometric data}

In \cite{V00}, Vorobjev considered
 what he called {\em a geometric data}.
 Here we will use the same terminology for a slightly different
 situation:

\medskip
 \noindent{\bf Definition.} Let
 $p : E \rightarrow P$ be a vector bundle and let
 ${\rm Vert}= {\rm ker} dp \subset TE$. A geometric data on
 the vector bundle $(E, p, P)$ consists of

\smallskip
$\bullet$ a connection $ \gamma: TE \rightarrow Vert,$

\smallskip
$\bullet$ a vertical bivector field  $\Pi^V,$

\smallskip
$\bullet$ and a 2-form $\mathbb F \in \Omega^2(P) \otimes  \cinf(E)$

\smallskip
\noindent such that
$$[\Pi^V, \Pi^V]=0; \leqno (i)$$ 
$$ [hor(u), \Pi^V]=0, \ \forall u \in \chi(P); \leqno(ii)$$
$$\partial_{\gamma} \mathbb F=0; \leqno (iii)$$
$${\rm Curv}_{\gamma}(u,v)= (\Pi^V)^\sharp(d \mathbb F (u,v)),  \
 \forall u,v \in \chi(P). \leqno(iv)$$

\noindent 
Unlike in \cite{V00}, we include the conditions
 $(i)$-$(iv)$ in the definition of a geometric data since we will
  consider only triples $(\gamma, \Pi^V,  \mathbb F)$
  satisfying those conditions. In fact,
 the main difference between Vorobjev's definition of a geometric data
  and the above one 
  is  that the 2-form 
 $\mathbb F$ is not necessarily nondegenerate. 
 Now, let us explain the above notations.
  Here $\gamma$ is an Ehresmann connection:
  at each point $e \in E $, $ \gamma_e: T_eE \rightarrow Vert_e$
 is a projection map. So $Hor:=$ker$\gamma$ gives an horizontal distribution. We have the splitting
 $$T_{e}E= Hor_{e} \oplus Vert_{e}, \quad \forall
  e \in E.$$ Consequently, for every vector field $u$ on the
 base manifold $P$, there is an horizontal vector field
 $hor(u)$  (tangent to $Hor$) which is obtained by lifting
 $u$. A 2-vector is ``vertical'' if it is a section of $\Lambda ^2Vert.$ The curvature of $\gamma$
  is given by
$${\rm Curv}_{\gamma}(u,v)= [hor(u), hor(v)]-
 hor[u,v], \ \forall u, v \in \chi(P).$$
 \noindent The operator $\partial_{\gamma}:
 \Omega^k(P) \otimes  \cinf(E) \rightarrow \Omega^{k+1}(P)
 \otimes  \cinf(E)$ is defined by
\begin{eqnarray*}
\partial_{\gamma} \mathbb G(u_0, \dots , u_{k})&=&
\sum_{i=0}^k (-1)^i \call_{hor(u)}(
 \mathbb G(u_0, \dots , \widehat{u_i} , \dots, u_{k}))\cr
 & & +\sum_{i<j} (-1)^{i+j}
\mathbb G([u_i, u_j], u_0, \dots ,\widehat{u_i} , \dots, \widehat{u_j} , \dots,  u_{k}).
\end{eqnarray*}

\noindent We have the following theorem
\begin{thm}
\label{geometric data}
 Fix a tubular neighborhood of a submanifold $P$ of a manifold $M$,
 it defines a vector bundle structure $p:E\rightarrow P$
 on an open neighborhood $E$ of $P$ ($P$ is identified to the zero section).
  Any Dirac structure on $M$ which has $P$ as a
 pre-symplectic leaf
 determines a geometric data on $E,$  up to a shrinking. 
\end{thm}

\noindent To prove Theorem \ref{geometric data}, we need to establish a couple
 of lemmas. We will first introduce  some notations.
Consider a point $m_0 \in P$ and a neighborhood
 $U$ of $m_0$ in $E$ with coordinates  $(x_i,\ y_j)$  as in Theorem
\ref{normal form} but with the additional condition that $x=$constant are the fibers of $p:E\rightarrow P.$
 Then, $Vert$ is generated by the vector fields $\frac{\partial}{\partial y_j}.$  Using the notations of
 Section 3, we have the following lemmas:

\begin{lemma}
\label{prep 3}
 For any $i, j, k \in \{1 , \dots , r\}$,
 we have $$ X_i \cdot \alpha_{jk}+
X_j \cdot \alpha_{ki}+X_k \cdot \alpha_{ij}=0.$$
\end{lemma}

\noindent {\it Proof:} We have the Courant bracket
 $$[ {\cal H}_i, {\cal H}_j]=
 ([X_i, X_j], \ \call_{X_i} \alpha_j - i_{X_j}d\alpha_i).$$
\noindent Since $L$ is isotropic, we get
 $$\langle  [ {\cal H}_i, {\cal H}_j], \
  {\cal H}_k \rangle=0.$$
 \noindent This gives
$$0= \alpha_k([X_i, X_j])+
 X_k \cdot(\alpha_j(X_i))+
(i_{X_i}d\alpha_j - i_{X_j}d\alpha_i)(X_k).$$
\noindent But $\alpha_k([X_i, X_j])=0$, for all
 $i, j, k \in \{1 , \dots , r\}$.
 Moreover,

$$d\alpha_j(X_i, X_k) -d\alpha_i(X_j, X_k)=
  X_i \cdot \alpha_{jk}
 -X_k \cdot \alpha_{ji}
 -X_j \cdot \alpha_{ik}
 +X_k \cdot \alpha_{ij}.$$
\noindent There follows
 $$0= \langle  [ {\cal H}_i, {\cal H}_j], \
  {\cal H}_k \rangle=  X_i \cdot \alpha_{jk}
 - X_j \cdot \alpha_{ik}
 +X_k \cdot \alpha_{ij}.$$
\noindent This completes the proof of the lemma.

\hfill \qed

\begin{lemma}
\label{prep 4}
 For any $i, j \in \{1 , \dots , r\}$,
 we have $$[ X_i, \  X_j]=  (\Pi^{V})^{\sharp}
  d\alpha_{ij}.$$
\end{lemma}

\noindent {\it Proof:}
Since $L$ is isotropic, we have
\begin{eqnarray*}
0&=&2 \langle  [ {\cal H}_i, {\cal H}_j], \
  {\cal V}_k \rangle \cr
 &=& dy_k([X_i, X_j])+
 Z_k \cdot \alpha_{ji}+
 (i_{X_i}d\alpha_j - i_{X_j}d\alpha_i) (Z_k).
\end{eqnarray*}
\noindent But
 $$d\alpha_j(X_i, Z_k)-
 d\alpha_i(X_j,  Z_k)=
  2 \sum_{\ell =1}^s Z_k^\ell {\partial \alpha_{ij} \over \partial
 y_{\ell}}= 2  Z_k \cdot\alpha_{ij}.$$

\noindent It follows that
$$2 \langle  [ {\cal H}_i, {\cal H}_j], \
  {\cal V}_k \rangle = dy_k([X_i, X_j])+ Z_k \cdot \alpha_{ji}+
 2  Z_k \cdot\alpha_{ij}= dy_k([X_i, X_j])+  Z_k \cdot\alpha_{ij}
.$$ \noindent We conclude that $$dy_k([X_i, X_j])= - \Pi^{V}(dy_k, d\alpha_{ij}).$$ There follows the lemma.

\hfill \qed

\medskip

\noindent {\it Proof of Theorem \ref{geometric data}:} 
 We first construct $\gamma $ or,
equivalently, the horizontal subbundle $Hor$. Define 
$$Hor_e = pr_1(L_e\cap (T_eM\oplus Vert_e^\circ)),$$
\noindent where $Vert_e^\circ$ is the annihilator of $Vert_e$.
In local coordinates as above, $Vert$ is generated by the vector fields $X_i$ and we have
$$hor(\frac{\partial}{\partial x_i})=X_i.$$
\noindent By definition, $\Pi^V$ is the bivector
 field given by putting together the Poisson 2-vectors we have on each fiber of
$E.$  Its local expressions are given by
$$\Pi^V(dy_i, dy_j)= - Z_j(x,y) \cdot y_i.$$
\noindent We define $\mathbb F$ by the formula
$$\mathbb F (u,v)=\Omega_L(hor(u),hor(v)).$$ Locally, we have the 
 components
$$\mathbb F \Big({\partial \over \partial x_i},
 {\partial \over \partial x_j}\Big)
= \alpha_i(X_j)= \alpha_{ij}= - \alpha_{ji}.$$

\noindent Notice that we may have to shrink $E$ in order  $\gamma$ and $\Pi^V$ be well defined.
Now, we have to show
 that Properties $(i)$-$(iv)$ of the definition of a geometric data hold. 
 But, it is sufficient to work in local coordinates. 
 Moreover, these properties are exactly  equivalent to
 the fact that $\Pi^V$ is  Poisson,
 Lemmas \ref{prep 2}, \ref{prep 3}, and \ref{prep 4}, respectively.
 Therefore,  any
 Dirac structure on $M$ having $P$ as pre-symplectic leaf 
 determines a geometric
data on the vector bundle $p: E \rightarrow P$
  (up to a shrinking of the total space $E$). 

 \hfill \qed

 Conversely, we have  the following theorem:
\begin{thm}
\label{Part 2}
 Any geometric data $(\gamma, \Pi^V, \mathbb F)$
  on a vector bundle $p: E \rightarrow P$ induces a
 Dirac structure on the total space $E.$
\end{thm}

We will divide the proof of Theorem \ref{Part 2} into lemmas.
 Suppose that $(\gamma, \Pi^V, \mathbb F)$ is
 a geometric data on a vector bundle $ (E, p, P)$.
 From now on, if $u$ is a vector field on $P$, we will simply denote by
  $\overline{u}$ the horizontal lift of $u$ (instead of $hor(u)$).
 Define 
 $$L^H_x={\rm Span}\Big\{ (\overline{u}, \alpha_u)_x \ | \ u \in \chi(P), \
 \alpha_{u | Vert}=0, \  \alpha_u(\overline{ v})= \mathbb F(u,v) \Big\},$$
$$L^V_x={\rm Span}\Big\{\Big((\Pi^V)^\sharp \beta,  \ \beta \Big)_x \ | \ 
 \beta_{|Hor}=0\Big\}.$$
Clearly, the subbundles $L^H$ and $L^V$  of $TM \oplus T^*M$ are isotropic 
 with respect to $\langle \cdot , \cdot \rangle $.
Moreover, we have the following lemma:

\begin{lemma}
\label{prep 5}
 Both spaces  $\Gamma(L^H)$ and $\Gamma(L^V)$ 
 are closed under the Courant bracket. 
\end{lemma}

\medskip

\noindent{\it Proof:} Let 
${\cal V}_u=(\overline{ u}, \alpha_u)$, ${\cal V}_v=(\overline{ v}, \alpha_v)$, and
 ${\cal V}_w=(\overline{ w}, \alpha_w)$ be elements of $\Gamma(L^H)$.
We have
\begin{eqnarray*}
2\langle [{\cal V}_u, {\cal V}_v], \ {\cal V}_w \rangle&=
& 2\langle ([\overline{ u}, \overline{ v}], \
 \call_{\overline{ u}}\alpha_v-i_{\overline{ v}}d \alpha_u), \ 
{\cal V}_w \rangle\cr
 &=& \Big(\alpha_w([\overline{ u}, \overline{ v}])+\call_{\overline{ u}}(\alpha_v(\overline{ w})) \Big)+c.p.,
\end{eqnarray*}
\noindent where the symbol $c.p.$ stands for the two other terms
 obtained by cyclic permutation of the indexes.
 Using the definition of $\alpha_v$ and the
  fact that  $$[\overline{ u}, \overline{ v}]= \overline{[u,v]}
 + (\Pi^V)^\sharp(d \mathbb F(u,v)),$$
\noindent  we obtain
$$ 2\langle [{\cal V}_u, {\cal V}_v], \ {\cal V}_w \rangle=
  \Big(\mathbb F(w, [u,v])+ \call_{\overline{ u}}(\mathbb F(v,w)) \Big) +c.p.
= \partial_{\gamma} \mathbb F(u,v,w)=0.$$

Similarly, one can show that the closedness of the space 
 $\Gamma(L^V)$ under the Courant 
bracket follows from the fact that $\Pi^V$ is a Poisson bivector field.

\hfill \qed

\begin{lemma}
\label{prep 6}
For any $u, v \in \chi(M)$,  $\beta \in (Hor)^\circ$, we have
 $$d\alpha_v(\overline{ u}, \  (\Pi^V)^\sharp \beta)= - \Pi^V(d\mathbb F(u,v), \beta).$$
\end{lemma}

\medskip

\noindent The proof of this lemma is straightforward. It is left to
 the reader.
\hfill \qed

\begin{lemma}
\label{prep 7}
For any ${\cal H}_1, {\cal H}_2 \in \Gamma(L^H)$,  and for any 
${\cal V} \in \Gamma(L^V)$, we have
$$\langle [{\cal H}_1, \ {\cal H}_2], \ {\cal V}
 \rangle=0.$$
\end{lemma}

 \noindent{\it Proof:} 
Let $${\cal H}_1=(\overline{ u}, \alpha_u), \quad 
{\cal H}_2=(\overline{ v}, \alpha_v) \quad {\rm and } \quad
 {\cal V}=((\Pi^V)^\sharp \beta, \beta).$$ 
\noindent Then,
\begin{eqnarray*}
2\langle [{\cal H}_1, \ {\cal H}_2], \ {\cal V} \rangle
&=& \beta\big([\overline{ u}, \overline{ v}]\big)
 +\Pi^V\big(\beta, \ d \mathbb F(u,v))\big)\cr
&& + d \alpha_v\big(\overline{ u}, \ (\Pi^V)^\sharp \beta \big)
 -d \alpha_u\big(\overline{ v}, \ (\Pi^V)^\sharp \beta \big) .
\end{eqnarray*}

\noindent Using Lemma \ref{prep 6}, we get
$$2\langle [{\cal H}_1, \ {\cal H}_2], \ {\cal V} \rangle=
 \beta([\overline{ u}, \overline{ v}])-\Pi^V(d\mathbb F(u,v), \beta)=0.$$

\hfill \qed

\begin{lemma}
\label{prep 8}
For any $u \in \chi(P)$, and for any
 1-forms $\beta_1, \beta_2$ on $M$ such that
 $\beta_{i_{|Hor}}=0$,  $i=1,2$, we have
$$d \beta_1 \Big((\Pi^V )^\sharp \beta_2, \ \overline{ u}\Big)=
  \call_{\overline{ u}}\Big((\Pi^V)^\sharp\beta_1, \
 \beta_2\Big)-\Pi^V\Big(\beta_1, \ \call_{\overline{ u}} \beta_2\Big).$$
\end{lemma}

\noindent {\it Proof:} Indeed, we have
$$d \beta_1\Big((\Pi^V)^\sharp \beta_2, \ \overline{ u}\Big)= 
- \call_{\overline{ u}}\Big(\Pi^V(\beta_1, \beta_2) \Big)+ 
 \beta_1\Big([\overline{ u}, \ (\Pi^V)^\sharp \beta_2]\Big).$$
Using the fact that $[\overline{ u}, \Pi^V]=0$, we get the formula of this lemma.

\hfill \qed

\begin{lemma}
\label{prep 9}
For any ${\cal V}_1, {\cal V}_2 \in \Gamma(L^V)$,  and for any 
${\cal H} \in \Gamma(L^H)$, we have
$$\langle [{\cal V}_1, \ {\cal V}_2], \ {\cal H}
 \rangle=0.$$
\end{lemma}
\noindent {\it Proof:}
Let ${\cal V}_i=((\Pi^V)^\sharp \beta_i, \ \beta_i)$, $i=1,2$, and ${\cal H}=(\overline{ u}, \alpha_u)$. By definition,

\begin{eqnarray*}
2\langle [{\cal V}_1, \ {\cal V}_2], \ {\cal H}
 \rangle&=& \alpha_u\Big([(\Pi^V)^\sharp \beta_1, \
 (\Pi^V)^\sharp \beta_2]\Big)
+ \call_{\overline{ u}}\Big(\Pi^V(\beta_1, \beta_2)\Big)\cr
 & & +
d\beta_2\Big((\Pi^V)^\sharp \beta_1, \ \overline{ u}\Big)-
d\beta_1\Big((\Pi^V)^\sharp \beta_2, \ \overline{ u}\Big).
\end{eqnarray*}

\noindent Now,
 using Lemma \ref{prep 8} and the fact that $\alpha_{u_{|Vert}}=0$,
we obtain $$2\langle [{\cal V}_1, \ {\cal V}_2], \ {\cal H}
 \rangle= - \call_{\overline{ u}}\Big(\Pi^V(\beta_1, \beta_2)\Big) +
\Pi^V(\beta_1, \ \call_{\overline{ u}}\beta_2)+
  \Pi^V(\call_{\overline{ u}}\beta_1, \ \beta_2) =0.$$

\hfill \qed

\begin{lemma}
\label{prep 10}
For any ${\cal V}_1, {\cal V}_2 \in \Gamma(L^V)$,  and for any 
${\cal H}_1, {\cal H}_2 \in \Gamma(L^H)$, we have
$$\langle [{\cal V}_1, \ {\cal H}_1], \ {\cal V}_2
 \rangle=0 \quad {\rm and} \quad 
\langle [{\cal V}_1, \ {\cal H}_1], \ {\cal H}_2
 \rangle=0.$$
\end{lemma}

\noindent {\it Proof:} Let
${\cal H}_i=(\overline{u}_i, \ \alpha_{u_i})$ 
 and ${\cal V}_i=((\Pi^V)^\sharp \beta_i, \ \beta_i)$, for $i=1,2$.
On the one hand, 
$$
2 \langle [{\cal H}_1, {\cal V}_1 ] , \ {\cal H}_2
\rangle= \alpha_{u_2}\Big([\overline{u}_1, \ (\Pi^V)^\sharp \beta_1]
 \Big)
+ d \beta_1(\overline{u}_1, \overline{u}_2)- d \alpha_{u_1} (\Pi^V(\beta_1), \overline{u}_2).$$

\noindent Since $[\overline{u}_1, \ \Pi^V]=0$, we get 
$$\alpha_{u_2}\Big([\overline{u}_1, \ (\Pi^V)^\sharp \beta_1]\Big)=
 \alpha_{u_2}\Big((\Pi^V)^\sharp \call_{\overline{u}_1}\beta_1 \Big)=0.$$
\noindent There follows 
\begin{eqnarray*}
2 \langle [{\cal H}_1, {\cal V}_1 ] , \ {\cal H}_2
\rangle
& =&  d \beta_1(\overline{u}_1, \overline{u}_2)- d \alpha_{u_1} (\Pi^V(\beta_1), \overline{u}_2) \cr 
&=& - \beta_1([\overline{u}_1, \overline{u}_2]) + \Pi^V(d \mathbb F(u_1, u_2), \beta_1) \ \mbox{by Lemma
\ref{prep 6}}\cr
&=&0 
\end{eqnarray*}
\noindent  since
 $[(\overline{u}_1, \overline{u}_2]= 
\overline{[u_1, u_2]} + (\Pi^V)^\sharp (d \mathbb F(u_1, u_2)). $
 On the other hand, 
\begin{eqnarray*}
2 \langle [{\cal H}_1, {\cal V}_1 ] , \ {\cal V}_2
\rangle&=& \beta_2\Big([\overline{u}_1, \ (\Pi^V)^\sharp \beta_1]\Big)
+ d \beta_1\Big(\overline{u}_1, \  (\Pi^V)^\sharp \beta_2\Big)\cr
&=& \Pi^V\Big(\call_{\overline{u}_1} \beta_1, \ \beta_2\Big) +
 \call_{\overline{u}_1} \Big(\Pi^V(\beta_2, \ \beta_1) \Big)
- \Pi^V\Big(\call_{\overline{u}_1} \beta_2, \ \beta_1\Big)\cr
 &=&0.
\end{eqnarray*} 

\noindent  This completes the proof of the lemma.

\hfill \qed
 
\medskip

\noindent   {\it Proof of Theorem \ref{Part 2}:}
Let $L$ be the vector bundle over $E$ whose fibre at $e \in E$
 is $L_e=L^H_e + L^V_e$. It follows
 immediately from Lemmas \ref{prep 5}, 
\ref{prep 7}, \ref{prep 9}, and \ref{prep 10} that $L$ is a Dirac
 structure on $E$.

\hfill \qed

\begin{cor}
 Any geometric data 
 $( \gamma, \Pi^V, \mathbb F)$ on a vector bundle
  $( E, p, P)$ whose  associated 2-form $\mathbb F$
is nondegenerate determines a Poisson structure on the total space $E$.
\end{cor}

\noindent{\it Proof:} Suppose
 $( \gamma, \Pi^V, \mathbb F)$ is a geometric data on
 $E$ such that $ \mathbb F$ is nondegenerate.
 Then, on each leaf of the foliation associated to the Dirac
 structure $L$ obtained from Theorem \ref{Part 2},
  the pre-symplectic  2-form  is nondegenerate.
 But, we known that all the leaves of a Dirac manifold are symplectic
 if and only if the associated
 Dirac structure is the graph of a Poisson bivector field.
 Hence, one gets the corollary.

\hfill \qed

\noindent{\bf Remark.} This
 corollary was given in \cite{V00} without proof.
 It has been established
 in \cite{Brahic} and \cite{Davis-W} by  methods
  which are different from the one used here.

\end{document}